\documentclass[preprint,prb,amsfonts,amsmath,amssymb]{revtex4}
\usepackage{bbding}
\usepackage{pifont}
\usepackage{graphicx}
\usepackage{amssymb}
\usepackage{txfonts}
\usepackage{dcolumn}
\usepackage{amsmath}
\usepackage{bm}

\usepackage{colordvi}

\begin{document}
\title{Several special cases of a square problem}
\author{Yang Ji}
\email{jiyang@semi.ac.cn}
\affiliation{1 State Key Laboratory for Superlattices and Microstructures, Institute of Semiconductors, Chinese Academy of Sciences, Beijing 100083, China \\
2 College of Materials Science and Opto-Electronic Technology, University of Chinese Academy of Sciences, Beijing 100049, China}


\begin{abstract}
Here is a square problem: in a unit square, is there a point with four rational distances to the vertices? A probability argument suggests a negative answer. This paper proves several special cases of the square problem: if the point sits on the diagonals, the midlines or the edges of the square, or the side-length of the square is $\mathrm{n}$ times the distance from the point to one side (both $n$ and $\left(n^{2}+4\right)$ are prime numbers), the distances from this point to the four vertices can not be all rational. However, this paper does not prove a more general situation. The proof here can be extended to the whole plane, instead of being limited to the interior of the square.

\end{abstract}

\maketitle

45 years ago, C.W. Dodge asked "a square problem " in the \textit{Mathematical Magazine}\cite{Dodge1976}: is there a point in a unit square which has all rational distances to the four vertices? He had no answer. 10 years later, in the same journal a comment\cite{MathMag1986}  said that John P. Robertson proved a special case:  there is no such a square if two of the distances from the point to the sides are equal. However, no details of the proof was given. A book published in 2005, \textit{Research Problems in discrete geometry} \cite{Brass2005}, still regarded the square problem as unproved.

This paper proves several special cases of the square problem: if the point is located on diagonal, midline or edge, or the side length of the square is  times the distance from the point to one side ( both $n$ and $\left(n^{2}+4\right)$  are prime numbers), this point cannot have four rational distances to the vertices of the unit square.

\begin{figure}[thb]
\includegraphics[width=0.55\columnwidth]{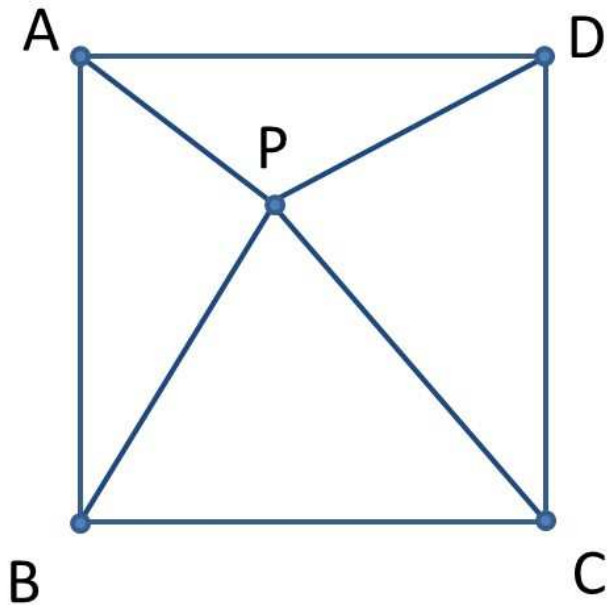}
\caption[Figure. 1] { A square problem. In a unit square, is there a point P with all rational distances to the four vertices? That is, the lengths of PA, PB, PC and PD are all rational numbers.}
\end{figure}

\section{Preliminary knowledge}

Suppose there is such a point $P$ in a unit square. Clearly, the coordinates of the point P must be rational. Multiplying with a suitable integer, we can get a square with side-width of an integer number, $z$. Its four vertices are $(0,0),(0, z),(z, z)$, and $(z, 0)$, respectively. The point $\mathrm{P}(x, y)$ has four integer distances to the sides, and four integer distances to the vertices. Assume $z$ is the smallest number with such
properties. Then, integers $x, y$ and $z$ have no common divisor.

Pythagorean triples $(p, q, r)$ can be expressed as $p=2 s t, q=s^{2}-t^{2}$ $r=s^{2}+t^{2}$, where $s$ and $t$ are odd and even numbers, respectively. Clearly, $x, y$ and $z$ should be combinations of some Pythagorean triples.

Because the sum of two squared odds cannot be a square, $z$ and one of $x$ and $y$ must be even. Assume it is $x$, then, $y$ must be odd. So, P is not on the diagonal, as also seen from the fact that the diagonal length is $\sqrt{2}$, which cannot be the sum of two rational numbers.

If the point is on the midline, the side-length of the square must have a form of 4st instead of $2\left(s^{2}-t^{2}\right)$. Else, the remainder of the side-length (moudulo 4) is 2, while the remainder of the sum of two $2 s t$ numbers is 0. A square cannot be like this.

If the side-length of the square is $n$ times the distance from that point to a certain edge $(n$ is an odd prime number $)$, this distance must be of the form $2 s t$, because $n-1$ is even, and an odd number cannot be even times of another odd number. Furthermore, the distances should be like $x=2 a b c d, \quad y=a^{2} c^{2}-b^{2} d^{2}$, and $z-y=b^{2} c^{2}-a^{2} d^{2}$, where $a$ is even, while $b, c$ and $d$ are odd. Thus, $z=\left(a^{2}+b^{2}\right)\left(c^{2}-d^{2}\right)$.

Consider the square problem with a probability argument. The numbers $x, y$ and $z$ are of the order of $a^{4}$, so are $z-x$ and $z-y$. The squared sum of two such numbers has a probability of to $a^{-4}$ be a square. Two such events have to occur simultaneously to satisfy the square problem , and the odd is thus $a^{-8}$. For a certain $a$, there are about $a^{3}$ kinds of possible $x$. Integration of $a$ leads to the total odd $\int_{a_{0}}^{\infty} a^{3} \cdot a^{-8} \mathrm{~d} a \sim a_{0}^{-4}$, where $a_{0}$
is starting point of the integration. A computer check can find that $a_{0}$ is bigger than 10 , or even 100 . So, it is a very small odd for a unit square to have a point with four rational distances to the vertices. However, this is not a proof.

\section{Proofs of several special cases}
It has been proved that there is no such a point on the diagonol. It is the same for the midline and the edge of the square, or when the side-length of the square is $n$ times the distance from the point to one side (both $n$ and $\left(n^{2}+4\right)$ are prime numbers), as the following theorems show.

\vspace{5 mm}

\textbf{Theorem 1.} On the edge of a unit square, no point has four rational distances to the vertices.

\textbf{Proof:} We use Fermat's method of infinite descent.

If there is such a point, then $x=2 a b c d=z=\left(a^{2}+b^{2}\right)\left(c^{2}-d^{2}\right)$, then, there are integers $a, b, c$ and $d, \quad$ which satisfy $2 a b c d=\left(a^{2}+\right.$ $\left.b^{2}\right)\left(c^{2}-d^{2}\right)$, where $a$ is even, while $b, c$ and $d$ are odd. Clearly, $a$ and $b$ are co-prime, so are $c$ and $d$.

Let $x=c / d, \quad k=a b /\left(a^{2}+b^{2}\right), \quad$ (please note that $x$ and $k$ are different variables, not the same as before), thus, the following eqation
$$
x^{2}-2 k x-1=0
$$
has a rational solution $x=c / d$, thus, its discriminant $4 k^{2}+4$ has to be a square number. It means the following equation has integer solutions,
$$
\begin{array}{c}
\left(a^{2}+b^{2}\right)^{2}+a^{2} b^{2}=e^{2} \quad\left({* }\right)
\end{array}
$$
where $e$ is an integer. Thus,
$$
\left(e-a^{2}-b^{2}\right)\left(e+a^{2}+b^{2}\right)=a^{2} b^{2}
$$
If $a_{1}$ is an odd factor of $a$, the two parts in the left cannot have a factor $a_{1}$ at the same time, else, we have
$$
\begin{array}{l}
\left(e-a^{2}-b^{2}\right)=k_{1} a_{1} \\
\left(e+a^{2}+b^{2}\right)=k_{2} a_{1}
\end{array}
$$
Say,
$$
e=\left(k_{1}+k_{2}\right) a_{1} / 2
$$
that is, $e$ has a factor $a_{1}$, thus, so does $b$, which means that both $a$ and
$b$ have a common factor $a_{1}$, contradictory to the requirement that they
are co-prime. Similarly, the odd factors of $b$ can only show in one of the
two parts in the left of the equation. Thus,
$$
\begin{array}{l}
\left(e-a^{2}-b^{2}\right)=k_{1} a_{1}^{2} b_{1}^{2} \\
\left(e+a^{2}+b^{2}\right)=k_{2} a_{2}^{2} b_{2}^{2}
\end{array}
$$
Where $a_{1}$ and $a_{2}\left(b_{1}\right.$ and $\left.b_{2}\right)$ are true factors of $a(b)$.

Because $a$ is even and $b$ is odd, the difference between $\left(e-a^{2}-b^{2}\right)$ and $\left(e+a^{2}+b^{2}\right)$ is $8 k+2$ (where $k$ is an integer), only one of them has a factor 2 and not 4, that is,
$$
\begin{array}{c}
\left(e-a^{2}-b^{2}\right)=2 a_{1}^{2} b_{1}^{2} \\
\left(e+a^{2}+b^{2}\right)=2^{2 k-1} a_{2}^{2} b_{2}^{2}
\end{array}
$$
where we assume $a=2^{k} a_{1} a_{2}$, both $a_{1}$ and $a_{2}$  are odd number, but not necessary to be prime. Thus,
$$
e=a_{1}^{2} b_{1}^{2}+2^{2 k-2} a_{2}^{2} b_{2}^{2}
$$
take this into Eq. $(*)$, we have
$$
\left(\left(2^{k} a_{1} a_{2}\right)^{2}+\left(b_{1} b_{2}\right)^{2}\right)^{2}+\left(2^{k} a_{1} a_{2} b_{1} b_{2}\right)^{2}=\left(\left(a_{1} b_{1}\right)^{2}+\left(2^{k-1} a_{2} b_{2}\right)^{2}\right)^{2}
$$
with the formula of the Pythagorean triples,
$$
\left(\left(a_{1} b_{1}\right)^{2}-\left(2^{k-1} a_{2} b_{2}\right)^{2}\right)^{2}+\left(2^{k} a_{1} a_{2} b_{1} b_{2}\right)^{2}=\left(\left(a_{1} b_{1}\right)^{2}+\left(2^{k-1} a_{2} b_{2}\right)^{2}\right)^{2}
$$
we have,
$$
\left(2^{k} a_{1} a_{2}\right)^{2}+\left(b_{1} b_{2}\right)^{2}=\left(a_{1} b_{1}\right)^{2}-\left(2^{k-1} a_{2} b_{2}\right)^{2}
$$
Change the symbols, $2^{k-1} a_{2} \rightarrow a, a_{1} \rightarrow c, \quad b_{1} \rightarrow b, \quad b_{2} \rightarrow d$, we have
$a^{2}\left(4 c^{2}+d^{2}\right)=b^{2}\left(c^{2}-d^{2}\right)$
that is,
$$
\frac{b^{2}}{a^{2}}=\frac{4 c^{2}+d^{2}}{c^{2}-d^{2}}
$$
so,
$$
\begin{array}{c}
k a^{2}=c^{2}-d^{2} \\
k b^{2}=4 c^{2}+d^{2}
\end{array}
$$
thus,
$$
5 c^{2}=k\left(a^{2}+b^{2}\right)
$$
$$
5 d^{2}=k\left(b^{2}-4 a^{2}\right)
$$
from which we have $k=5 .$ Else, $c$ and $d$ have co-factors, leading to a contradiction. So,
$$
\begin{array}{c}
a^{2}+b^{2}=c^{2} \quad\left({**}\right) \\
d^{2}+4 a^{2}=b^{2} \quad\left({***}\right)
\end{array}
$$
Let $d=m^{2}-n^{2}, a=m n, b=m^{2}+n^{2}$, Eq. $(***)$ holds, and Eq. $(**)$ turns into
$$
\left(m^{2}+n^{2}\right)^{2}+(m n)^{2}=c^{2}
$$
Again change the variables, $m \rightarrow a, n \rightarrow b, \quad c \rightarrow e$, we go back to the
starting equation,
$$
\begin{array}{c}
\left(a^{2}+b^{2}\right)^{2}+a^{2} b^{2}=e^{2} \quad\left({* }\right)
\end{array}
$$
This makes Fermat's infinite descent. So, the equation has no integer solution.

Q. E. D.

\vspace{5 mm}

\textbf{Theorem 2.} On the midline of a unit square, no point has four rational distances to the vertices.

\textbf{Proof: }Again we use Fermat's method of infinite descent.

Similar to Theorem 1, but now we will show the following equantion has no integer solution.
$$
\left(a^{2}+b^{2}\right)^{2}+(2 a b)^{2}=e^{2}
$$
that is,
$$
\left(\left(a_{1} a_{2}\right)^{2}+\left(b_{1} b_{2}\right)^{2}\right)^{2}+\left(2 a_{1} a_{2} b_{1} b_{2}\right)^{2}=e^{2}
$$
With similar procedures of Theorem 1 , we have to show
$$
\frac{b^{2}}{a^{2}}=\frac{c^{2}+d^{2}}{c^{2}-d^{2}}
$$
has no integer solution, where $a$ is even. Thus,
$$
\begin{array}{l}
2 c^{2}=k\left(a^{2}+b^{2}\right) \\
2 d^{2}=k\left(b^{2}-a^{2}\right)
\end{array}
$$
from which we have $k=2$. Else, $c$ and $d$ have co-factors, leading to a contradiction. So,
$$
\begin{array}{l}
a^{2}+b^{2}=c^{2} \\
d^{2}+a^{2}=b^{2}
\end{array}
$$
Let $b=m^{2}+n^{2}, \quad a=2 m n, \quad d=m^{2}-n^{2}$, the above two equations lead to
the following,
$$
\left(m^{2}+n^{2}\right)^{2}+(2 m n)^{2}=c^{2}
$$
Change the variables, $m \rightarrow a, n \rightarrow b, c \rightarrow e$, we go back to the starting equation,
$$
\left(a^{2}+b^{2}\right)^{2}+(2 a b)^{2}=e^{2}
$$
Again we get Fermat's infinite descent.

Q. E. D.

\vspace{5 mm}

\textbf{Theorem 3.} If the side-length is $n$ times the distance from a point to one side of the square ( both $n$ and $\left(n^{2}+4\right)$ are prime numbers), this point cannot have four rational distances to the vertices of the square.

\textbf{Proof: }Similar to Theorem 1 and 2 , but now we will show the following equation has no integer solution.
$$
\left(a^{2}+b^{2}\right)^{2}+(n a b)^{2}=e^{2}
$$
that is,
$$
\left(\left(a_{1} a_{2}\right)^{2}+\left(b_{1} b_{2}\right)^{2}\right)^{2}+\left(n a_{1} a_{2} b_{1} b_{2}\right)^{2}=e^{2}
$$
Again, we have to prove,
$$
\frac{b^{2}}{a^{2}}=\frac{4 c^{2}+n^{2} d^{2}}{c^{2}-d^{2}}
$$
has no integer solution, where $a$ is even. Thus,
$$
\begin{array}{l}
\left(n^{2}+4\right) c^{2}=k\left(n^{2} a^{2}+b^{2}\right) \\
\left(n^{2}+4\right) d^{2}=k\left(b^{2}-4 a^{2}\right)
\end{array}
$$
If $n^{2}+4$ is a prime number, we have $k=n^{2}+4$. Else, $c$ and $d$ have co-factors, leading to a contradiction. So,
$$
\begin{array}{c}
n^{2} a^{2}+b^{2}=c^{2} \\
b^{2}-4 a^{2}=d^{2}
\end{array}
$$
Let $b=s^{2}+t^{2}, a=s t, \quad d=s^{2}-t^{2}$, the above two equations lead to the following,
$$
\left(s^{2}+t^{2}\right)^{2}+(n s t)^{2}=c^{2}
$$
Change the variables, $s \rightarrow a, t \rightarrow b, c \rightarrow e$, we go back to the starting equation,
$$
\left(a^{2}+b^{2}\right)^{2}+(n a b)^{2}=e^{2}
$$
Again we get Fermat's infinite descent.

Q. E. D.

\vspace{5 mm}

\textbf{Remark 1: }The above theorems do not require the point P to be in the square. It can be applied to all points in the plane. According to the previous probability arguments, in the whole plane, there may be no point with four rational distances to the vertices of the unit square.

\textbf{Remark 2:} There are infinite points with three rational distances to the vertices of a unit square. Extend the shortest side of the Pythagorean triangle to the length of the other right-angle side, and then make a copy of the new triangle, these two will form the square that is needed. Shute and Yocom said in 1977,  there is a parameter formula that could give infinite points with three rational distances \cite{Shute1977}. However, they did not give details.

\section{Conclusion}

With Fermat's method of infinite descent, we proved several special cases of the square problem: if the point sits on the diagonals, the midlines or the edges of the square, or the side-length of the square is n times the distance from the point to one side (both  and  are prime numbers), the distances from this point to the four vertices can not be all rational. However, this paper does not prove a more general situation. The proof here can be extended to the whole plane, instead of being limited to the interior of the square.

\begin{acknowledgements}
 I thank ZhiYeShuXueJiaZaiMinJian and FengYunBianHuan for telling me the square problem and its references. I also thank ChinaXiv for comments on the paper.
 
 \vspace{5mm}
 
\textbf{ Note added in the revised version.} I thank Michel Marcus for telling me the following information.

1. Theorem 1  (for the boundary points) is not new. See, 95.01 The rational distance problem, ROY BARBARA, The Mathematical Gazette, Vol. 95, No. 532 (March 2011), pp. 59-61 (3 pages)
Published By: The Mathematical Association. https://www.jstor.org/stable/23248619.

2. For the parameter formula mentioned in Remark 2, see A215365 of the OEIS (https://oeis.org/A215365).

3. Theorem 3 needs both $p$ and $p^2+4$ to be prime number. For such prime pairs, see  A062324 and A045637 of the OEIS. The number of such prime pairs is infinite under the Bunyakovsky conjecture. - Charles R Greathouse IV, Jul 04 2011.

\end{acknowledgements}

\clearpage

\end{document}